\documentclass{birkjour}
\usepackage{mathtools}
 \newtheorem{thm}{Theorem}[section]
 \newtheorem{cor}[thm]{Corollary}
 \newtheorem{lem}[thm]{Lemma}
 
 \theoremstyle{definition}
 \newtheorem{defn}[thm]{Definition}
 \theoremstyle{remark}
 
 \newtheorem{ex}{Example}[section]
 \numberwithin{equation}{section}
\begin{document}

\title[$*$-Conformal $\eta$-Ricci soliton on Kenmotsu manifolds]
{$*$-Conformal $\eta$-Ricci soliton within the framework of Kenmotsu manifolds}

\author[S. Sarkar]{Sumanjit Sarkar}
\address{Department of Mathematics\\
Jadavpur University\\
Kolkata-700032, India.}
\email{imsumanjit@gmail.com}

\author[S. Dey]{Santu Dey}
\address{Department of Mathematics\\
Bidhan Chandra College\\
Asansol, Burdwan, West Bengal-713304, India.}
\email{santu.mathju@gmail.com}

%\author[A. Bhattacharyya]{Arindam Bhattacharyya}
%\address{Department of Mathematics\\
%Jadavpur University\\
%Kolkata-700032, India}
%\email{bhattachar1968@yahoo.co.in}


\begin{abstract}
The goal of our present paper is to deliberate $*$-conformal $\eta$-Ricci soliton within the framework of Kenmotsu manifolds. Here we have shown that a Kenmotsu metric as a $*$-conformal $\eta$-Ricci soliton is Einstein metric if the soliton vector field is contact. Further, we have evolved the characterization of the Kenmotsu manifold or the nature of the potential vector field when the manifold satisfies gradient almost $*$-conformal $\eta$-Ricci soliton. Next, we have contrived $*$-conformal $\eta$-Ricci soliton admitting $(\kappa,\mu)^\prime$-almost Kenmotsu manifold and proved that the manifold is Ricci flat and is locally isometric to $\mathbb{H}^{n+1}(-4)\times\mathbb{R}^n$. Finally we have constructed some examples to illustrate the existence of $*$-conformal $\eta$-Ricci soliton, gradient almost $*$-conformal $\eta$-Ricci soliton on Kenmotsu manifold and $(\kappa,\mu)^\prime$-almost Kenmotsu manifolds.
\end{abstract}

\maketitle

\section{\textbf{Introduction}}
In modern mathematics, the methods of contact geometry play an important role. Contact geometry has evolved from the mathematical formalism of classical mechanics. In 1969, S. Tanno\cite{Tanno} classified the connected almost contact metric manifolds whose automorphism groups have maximal dimensions as follows.\par
(i) Homogeneous normal contact Riemannian manifolds with constant $\phi$-holomorphic sectional curvature if $k(\xi,X)>0$;\par
(ii) Global Riemannian product of a line or a circle and a K\"{a}hlerian manifold with constant holomorphic sectional curvature if $k(\xi,X)=0$;\par
(iii) A warped product space $\mathbb{R}\times_fN$, where $\mathbb{R}$ is the real line and $N$ is a K\"{a}hlerian manifold, if $k(\xi,X)<0$;\\
where $k(\xi,X)$ denotes the sectional curvature of the plane section containing the characteristic vector field $\xi$ and an arbitrary vector field $X$.\par

In 1972, K. Kenmotsu in \cite{Kenmotsu} obtained some tensor equations to characterize the manifolds of the third class using the wraping function $f(t)=ce^t$ on the interval $J=(-\epsilon,\epsilon)$. Since then the manifolds of the third class were called Kenmotsu manifolds. Conversely, every point on a Kenmotsu manifold has a neighbourhood which is locally a warped product $J\times_fN$, where $f$ is given by the above mentioned relation.\par

A pseudo-Riemannian manifold $(M,g)$ admits a Ricci soliton which is a generalization of Einstein metric (i.e, $S=ag$ for some constant $a$) if there exists a smooth vector field $V$ and a constant $\lambda$ such that,
\begin{equation*}
  \frac{1}{2}\mathcal{L}_Vg+S+\lambda g=0,
\end{equation*}
where $\mathcal{L}_V$ denotes Lie derivative along the direction $V$ and $S$ denotes the Ricci curvature tensor of the manifold. The vector field $V$ is called potential vector field and $\lambda$
is called soliton constant.\par

The Ricci soliton is a self-similar solution of the Hamilton's Ricci flow \cite{Ham} which is defined by the equation $\frac{\partial g(t)}{\partial t}=-2S(g(t))$ with initial condition $g(0)=g$ where $g(t)$ is a one-parameter
family of metrices on $M$. The potential vector field $V$ and soliton constant $\lambda$ play  vital roles while determining the nature of the soliton. A soliton is said to be shrinking, steady or expanding according as $\lambda<0$,
$\lambda=0$ or $\lambda>0$. Now if $V$ is zero or Killing then the Ricci soliton reduces to Einstein manifold and the soliton is called trivial soliton.\par

If the potential vector field $V$ is the gradient of a smooth function $f$, denoted by $Df$ then the soliton equation reduces to,
\begin{equation*}
  Hessf+S+\lambda g=0,
\end{equation*}
where $Hessf$ is Hessian of $f$. Perelman\cite{Perelman} proved that a Ricci soliton on a compact manifold is a gradient Ricci soliton.\par

In 2005, A. E. Fischer \cite{Fischer} has introduced conformal Ricci flow which is a mere generalisation of the classical Ricci flow equation that modifies the unit volume constraint to a scalar curvature constraint. The conformal Ricci flow equation was given by,
\begin{eqnarray*}
  \frac{\partial g}{\partial t}+2(S+\frac{g}{n}) &=& -pg, \\
  r(g) &=& -1,
\end{eqnarray*}
where $r(g)$ is the scalar curvature of the manifold, $p$ is scalar non-dynamical field and $n$ is the dimension of the manifold. Corresponding to the conformal Ricci flow equation in 2015, N. Basu and A. Bhattacharyya \cite{Basu} introduced the notion of conformal Ricci soliton equation as a generalization of Ricci soliton equation given by,
\begin{equation*}
\mathcal{L}_Vg+2S+[2\lambda-(p+\frac{2}{n})]g=0.
\end{equation*}

In 2009, J. T. Cho and M. Kimura \cite{Kimura} introduced the concept of $\eta$-Ricci soliton which is another generalization of classical Ricci soliton  and is given by,
\begin{equation*}
  \mathcal{L}_\xi g+2S+2\lambda g+2\mu\eta\otimes\eta=0,
\end{equation*}
where $\mu$ is a real constant, $\eta$ is a 1-form defined as $\eta(X)=g(X,\xi)$ for any $X\in\chi(M)$. Clearly it can be noted that if $\mu=0$ then the $\eta$-Ricci soliton reduces to Ricci soliton.\par

Recently Md. D. Siddiqi \cite{Siddiqi} established the notion of conformal $\eta$-Ricci soliton which generalizes both conformal Ricci soliton and $\eta$-Ricci soliton. The equation for conformal $\eta$-Ricci soliton is given by,
\begin{equation*}
\mathcal{L}_\xi g+2S+[2\lambda-(p+\frac{2}{n})]g+2\mu\eta\otimes\eta=0.
\end{equation*}

In 2014, Kaimakamis and Panagiotidou \cite{Kai} modified the definition of Ricci soliton where they have used $*$-Ricci tensor $S^*$ which was introduced by Tachibana \cite{Tachi}, in place of Ricci tensor $S$. The $*$-Ricci tensor $S^*$ is defined by,
\begin{equation*}
  S^*(X,Y)=\frac{1}{2}(trace\{\phi\cdot R(X,\phi Y)\}),
\end{equation*}
for all vector fields $X$ and $Y$ on $M$, where $\phi$ is a $(1,1)-$tensor field. They have used the concept of $*$-Ricci soliton within the framework of real hypersurfaces of a complex space form.\par

In 2019, S. Roy et al.\cite{Soumendu} defined $*$-conformal $\eta$-Ricci soliton as,
\begin{equation*}
  \mathcal{L}_\xi g+2S^*+[2\lambda-(p+\frac{2}{n})]g+2\mu\eta\otimes\eta=0.
\end{equation*}

As per the authors knowledge, the results concerning $*$-conformal $\eta$-Ricci soliton were studied when the potential vector field $V$ is the characteristic vector field $\xi$. Motivated from this we generalize the definition by considering the potential vector field as arbitrary vector field $V$ and define as:
\begin{equation}\label{A1}
  \mathcal{L}_Vg+2S^*+[2\lambda-(p+\frac{2}{(2n+1)})]g+2\mu\eta\otimes\eta=0,
\end{equation}
where we considered the manifold as $(2n+1)$-dimensional. Now if we consider the potential vector field $V$ as the gradient of a smooth function $f$, then the $*$-conformal $\eta$-Ricci soliton equation can be rewritten as,
\begin{equation}\label{A2}
  Hess f+S^*+[\lambda-(\frac{p}{2}+\frac{1}{(2n+1)})]g+\mu\eta\otimes\eta=0.
\end{equation}

By gradient almost $*$-conformal $\eta$-Ricci soliton we mean gradient $*$-conformal $\eta$-Ricci soliton where we consider $\lambda$ as a smooth function.\par
It is worthy to mention that Sharma\cite{Sharma} first initiated the study of Ricci solitons in contact geometry. However, Ghosh\cite{Ghosh} is the first to consider 3-dimensional Kenmotsu metric as a Ricci soliton. After that Kenmotsu manifold is studied on many context of Ricci soliton by many authors like: C\v{a}lin and Crasmareanu \cite{Calin}, Ghosh \cite{Ghosh1}, Wang \cite{Wang} etc. In \cite{Naik}, authors have considered $*$-Ricci solitons and gradient almost $*$-Ricci solitons on Kenmotsu manifolds and obtained some beautiful results. They have proved that if a 3-Kenmotsu manifold admits a $*$-Ricci soliton, then the manifold is of constant negative curvature $-1$. They also get that if an $\eta$-Einstein Kenmotsu manifold of $dim>3$ admits a $*$-Ricci soliton then the manifold becomes Einstein. For the Kenmotsu manifold possessing a gradient almost $*$-Ricci soliton they come up with: either the manifold is Einstein or the potential vector field is pointwise collinear with the characteristic vector field on some open set of the manifold. Further, Dey et al. \cite{dey1} have studied $*$-$\eta$-Ricci soliton on contact geometry. Being motivated from the well acclaimed results we consider $*$-conformal $\eta$-Ricci soliton and gradient almost $*$-conformal $\eta$-Ricci soliton in the framework of Kenmotsu manifolds.\par

\subsection{\textbf{Physical motivation:}}
The notion of $*$-conformal $\eta$-Ricci soliton is replaced by conformal $\eta$-Ricci soliton as a kinematic solution in fluid space-time, whose profile develop a characterization of spaces of constant sectional curvature along with the locally symmetric spaces. The symmetries are the part of geometry and thus reveals the physics. The metric symmetries are important as they simplify solutions to many problems. Their main application in general relativity is that they classify solutions of Einstein field equations. One of these symmetries is $*$-conformal $\eta$-Ricci soliton associated with space-time geometry. Also geometric phenomenon of $*$-conformal $\eta$-Ricci solitons can evolve an aqueduct between a sectional curvature inheritance symmetry of space-time although metrics considered in the work are positivedefinite and class of conformal $\eta$-Ricci solitons.\par

The mathematical notion of a $*$-conformal $\eta$-Ricci soliton should not be confused with the notion of soliton solutions, which arise in several areas of methematical or theoretical physics and its applications. $*$-conformal $\eta$-Ricci soliton is important as it can help in understanding the concepts of energy or entropy in general realtivity. This property is the same as that of heat equation due to which an isolated system loses the heat for a thermal equilibrium. By \cite{Duggal},\cite{Woolgar} (briefly discussed in the above) we can think more about physical applications of $*$-conformal $\eta$-Ricci soliton.\\ \par

The paper is organized as follows: in section 2, the basic definitions and facts about contact metric manifolds, Kenmotsu manifolds and $(\kappa,\mu)^\prime$-almost Kenmotsu manifolds are given. In the later section, we consider Kenmotsu metric as $*$-conformal $\eta$-Ricci soliton and gradient almost $*$-conformal $\eta$-Ricci soliton and obtain some useful results. We also provide some examples to support our findings in that section. In section 4, we consider the metric of $(\kappa,\mu)^\prime$-almost Kenmotsu manifold to represent $*$-conformal $\eta$-Ricci soliton along with a special condition and obtained that the manifold is Ricci flat and is locally isometric to $\mathbb{H}^{n+1}(-4)\times\mathbb{R}^n$. We also furnish an example at the end of this section. The last section comes up with some physical motivations and conclusion.

\section{\textbf{Notes on contact metric manifolds}}
By \cite{Blair}, a differentiable manifold $M$ of dimension $(2n+1)$ is said to has an almost contact structure or $(\phi,\xi,\eta)$ structure if $M$ admits a $(1,1)$ tensor field $\phi$, a vector field $\xi$, a 1-form $\eta$ satisfying:
\begin{eqnarray}
  \phi ^2 &=& -I + \eta \otimes \xi, \\ \label{B1}
  \eta (\xi ) &=& 1, \label{B2}
\end{eqnarray}
where $I$ is the identity mapping. Generally, $\xi$ and $\eta$ are called \emph{characteristic vector field} or \emph{Reeb vector field} and \emph{almost contact 1-form} respectively.\par
A Riemannian metric $g$ is said to be \emph{compatible metric} if it satisfies:
\begin{equation}
  g( \phi X , \phi Y) = g(X,Y)- \eta (X) \eta (Y), \label{B3}
\end{equation}
for arbitrary vector fields $X$ and $Y$ on $M$. A manifold having almost contact structure along with compatible Riemannian metric is called \emph{almost contact metric manifold}.\par
In a almost contact metric manifold $(M,\phi,\xi,\eta,g)$ the following conditions are satisfied:
\begin{eqnarray}
  \phi \xi &=& 0, \label{B4}\\
  \eta \circ \phi &=& 0, \label{B5}\\
  g(X, \xi ) &=& \eta (X), \label{B6}\\
  g(\phi X,Y) &=& -g(X,\phi Y),\label{B7}
\end{eqnarray}
for arbitrary $X,Y \in \chi (M)$. The normality of an almost contact structure is equivalent with the vanishing of the tensor $N_\phi=[\phi,\phi]+2d\eta\otimes\xi$, where $[\phi,\phi]$ is the Nijenhuis tensor of $\phi$ (for more details we refer to \cite{Blair}).
\begin{defn}
On an almost contact metric manifold $M$, a vector field $X$ is said to be contact vector field if there exist a smooth function $f$ such that $\mathcal{L}_X\xi=f\xi$.
\end{defn}
\begin{defn}
On an almost contact metric manifold $M$, a vector field $X$ is said to be infinitesimal contact transformation if $\mathcal{L}_X\eta=f\eta$, for some function $f$. In particular, we call $X$ as a strict infinitesimal contact transformation if $\mathcal{L}_X\eta=0$.
\end{defn}
\medskip
\medskip

\subsection{\textbf{Kenmotsu manifold:}}
We define the fundamental 2-form $\Phi$ on an almost contact metric manifold $M$ by $\Phi(X,Y)=g(X,\phi Y)$ for arbitrary $X,Y\in\chi(M)$. We recall from \cite{Janssens}, an \emph{almost Kenmotsu manifold} is an almost contact metric manifold where $\eta$ is closed, i.e., $d\eta=0$ and $d\Phi=2\eta\wedge\Phi$. A normal almost Kenmotsu manifold is called \emph{Kenmotsu manifold}.\par
By \cite{Blair} if in a almost contact metric manifold $M$ the 1-form $\eta$ and the (1,1)-tensor field $\phi$ satisfy the following condition for arbitrary $X,Y \in \chi (M)$:
 \begin{equation}
   (\nabla_X \phi )Y = g(\phi X,Y)\xi-\eta(Y)\phi X,\label{B8}
 \end{equation}
 where $\nabla$ denotes the Riemannian connection of $g$, then the manifold $M$ is called a Kenmotsu manifold. It is easy to verify that the above mentioned relation is equivalent with the normality condition of the manifold.
 \par
 In Kenmotsu manifold of dimension $(2n+1)$ the following relationss hold:
 \begin{eqnarray}
   \nabla_X \xi &=& X- \eta(X)\xi,\label{B9}\\
   (\nabla_X \eta)Y &=& g(X,Y)-\eta(X)\eta(Y), \label{B10}\\
   R(X,Y)\xi &=& \eta(X)Y-\eta(Y)X, \label{B11}\\
   S(X,\xi) &=& -2n\eta (X), \label{B12}\\
   \mathcal{L}_\xi g(X,Y) &=& 2g(X,Y)-2\eta(X)\eta(Y),\label{B13}
 \end{eqnarray}
 for arbitrary $X,Y,Z,W\in \chi (M)$, where $\mathcal{L}$ is the Lie derivative operator, $R$ is Riemannian curvature tensor and $S$ is the Ricci tensor.\par
 An (2n+1)-dimensional Kenmotsu metric manifold is said to be $\eta$-Einstein Kenmotsu manifold if there exists two smooth functions $a$ and $b$ which satisfies the following relation,
\begin{equation}
  S(X,Y)=ag(X,Y)+b\eta(X)\eta(Y),\label{B14}
\end{equation}
for all $X,Y\in\chi(M)$. Clearly if $b=0$ then $\eta$-Einstein manifold reduces to Einstein manifold. Now considering $X=\xi$ in the last equation and using (\ref{B12}) we have, $a+b=-2n$. Contracting (\ref{B14}) over $X$ and $Y$ we get, $r=(2n+1)a+b$ where $r$ denotes the scalar curvature of the manifold. Solving these two we have, $a=(1+\frac{r}{2n})$ and $b=-(2n+1+\frac{r}{2n})$. Using these values we can rewrite (\ref{B14}) as,
\begin{equation}\label{B15}
  S(X,Y)=(1+\frac{r}{2n})g(X,Y)-(2n+1+\frac{r}{2n})\eta(X)\eta(Y).
\end{equation}

\subsection{($\kappa,\mu)^\prime$ almost Kenmotsu manifold:}
On an almost Kenmotsu manifold we consider two (1,1)-type tensor fields $h=\frac{1}{2}\mathcal{L}_\xi\phi$ and $h^\prime=h\circ\phi$ and an operator $\ell=R(.,\xi)\xi$, where $\mathcal{L}_\xi\phi$ is the Lie derivative of $\phi$ along the direction $\xi$. The tensor fields $h$ and $h^\prime$ plays vital role in almost Kenmotsu manifold. Both of them are symmetric and satisfy the following relations:
\begin{eqnarray}
  \nabla_X \xi &=& X-\eta(X)\xi+h^\prime X, \label{B18}\\
  h\xi&=&h^\prime\xi=0, \label{B16}\\
  h\phi &=& -\phi h, \label{B17}\\
  tr(h)&=&tr(h^\prime)=0,\nonumber
\end{eqnarray}
for any $X,Y\in\chi(M)$, where $\nabla$ is the Levi-Civita connection of the metric $g$. In addition the following curvature property is aslo satisfied:
\begin{equation}\label{B19}
  R(X,Y)\xi = \eta(X)(Y+h^\prime Y)-\eta(Y)(X+h^\prime X)+(\nabla_X h^\prime)Y-(\nabla_Y h^\prime)X,
\end{equation}
where $R$ is the Riemannian curvature tensor of $(M,g)$.\par
By $(\kappa,\mu)^\prime$-almost Kenmotsu manifold we mean almost Kenmotsu manifold where the characteristic vector field $\xi$ satisfies the $(\kappa,\mu)^\prime$-nullity distribution (for details see \cite{Dileo}), i.e.,
\begin{equation}\label{B20}
  R(X,Y)\xi=\kappa(\eta(Y)X-\eta(X)Y)+\mu(\eta(Y)h^\prime X-\eta(X)h^\prime Y),
\end{equation}
for any $X,Y\in\chi(M)$, where $\kappa$ and $\mu$ are real constants. On a $(\kappa,\mu)^\prime$-almost Kenmotsu manifold $M$ we have (see \cite{Dileo}),
\begin{eqnarray}
  h^{\prime2}(X)&=&-(\kappa+1)[X+\eta(X)\xi],\label{B21}\\
  h^2(X)&=&-(\kappa+1)[X+\eta(X)\xi],
\end{eqnarray}
for $X\in\chi(M)$. From previous relation it follows that $h^\prime=0$ if and only if $\kappa=-1$ and $h^\prime\neq0$ otherwise. Let $X\in Ker(\eta)$ be an eigenvector field of $h^\prime$ orthogonal to $\xi$ w.r.t. the eigenvalue $\alpha$. Then, from (\ref{B21}) we get $\alpha^2=-(\kappa+1)$ which implies $\kappa\leq -1$. Dileo and Pastore proved that on a $(\kappa,\mu)^\prime$-almost Kenmotsu manifold with $\kappa<-1$, we have $\mu=-2$ (Proposition 4.1 of \cite{Dileo}). Since the same symbol $\mu$ is used in the coefficient of $\eta\otimes\eta$ in the definition of $*$-conformal $\eta$-Ricci soliton and in $(\kappa,\mu)^\prime$-almost Kenmotsu manifold, so to reduce the complications in notations we use $(\kappa,-2)^\prime$-almost Kenmotsu manifold throughout this paper.\par
We recall some useful results on a $(2n+1)$ dimensional $(\kappa,-2)^\prime$-almost Kenmotsu manifold $M$ with $\kappa<-1$ as follows:
\begin{eqnarray}
  R(\xi,X)Y &=& \kappa(g(X,Y)\xi-\eta(Y)X)-2(g(h^\prime X,Y)\xi-\eta(Y)h^\prime X)~~ \label{B22}\\
  QX &=& -2nX+2n(\kappa+1)\eta(X)\xi-2nh^\prime(X),\label{B23} \\
  r &=& 2n(\kappa-2n),\label{B24}\\
  (\nabla_X\eta)Y &=& g(X,Y)-\eta(X)\eta(Y)+g(h^\prime X,Y), \label{B25}
\end{eqnarray}
where $X,Y\in\chi(M)$, $Q,r$ are the Ricci operator and scalar curvature of $M$ respectively.

\section{\textbf{$*$-conformal $\eta$-Ricci soliton on Kenmotsu manifold}}
In this section we consider that the metric $g$ of a $(2n+1)$-dimensional Kenmotsu manifold represents a $*$-conformal $\eta$-Ricci soliton and a gradient almost $*$-conformal $\eta$-Ricci soliton. We recall some important lemmas relevant to our results.

\begin{lem}\cite{Naik}
The Ricci operator $Q$ on a $(2n+1)$-dimensional Kenmotsu manifold satisfies
\begin{eqnarray}
  (\nabla_XQ)\xi &=& -QX-2nX, \label{C1}\\
  (\nabla_\xi Q)X &=& -2QX-4nX, \label{C2}
\end{eqnarray}
for arbitrary vector field $X$ on the manifold.
\end{lem}

\begin{lem}\cite{Naik}
The $*$-Ricci tensor $S^*$ on a $(2n+1)$-dimensional Kenmotsu manifold is given by,
\begin{equation}\label{C3}
  S^*(X,Y)=S(X,Y)+(2n-1)g(X,Y)+\eta(X)\eta(Y),
\end{equation}
for arbitrary vector fields $X$ and $Y$ on the manifold.
\end{lem}

\begin{thm}
Let $M^{(2n+1)}(\phi,\xi,\eta,g)$ be a Kenmotsu manifold. If the metric $g$ represents a $*$-conformal $\eta$-Ricci soliton and if the soliton vector field $V$ is contact, then V is strictly infinitesimal contact transformation and the manifold is Einstein.
\begin{proof}
  Since the metric $g$ of the Kenmotsu manifold represents a $*$-conformal $eta$-Ricci soliton so both of the equations (\ref{A1}) and (\ref{C3}) are satisfied. Combining these two we have,
  \begin{eqnarray}
    (\mathcal{L}_Vg)(X,Y)&=&-2S(X,Y)-(2\lambda-p-\frac{2}{(2n+1)}+4n-2)g(X,Y)\nonumber\\
    &&-2(\mu+1)\eta(X)\eta(Y). \label{C4}
  \end{eqnarray}
  Taking covariant derivative w.r.t. arbitrary vector field $Z$ and using (\ref{B10}) we have,
  \begin{eqnarray}
    (\nabla_Z\mathcal{L}_Vg)(X,Y) &=& -2(\nabla_ZS)(X,Y)-2(\mu+1)\{g(X,Z)\eta(Y) \nonumber\\
    && +g(Y,Z)\eta(X)-2\eta(X)\eta(Y)\eta(Z)\},\label{C5}
  \end{eqnarray}
      for all $X,Y,Z\in\chi(M)$. Again from Yano\cite{Yano} we have the following commutation formula,
    \begin{eqnarray}
      (\mathcal{L}_V\nabla_Zg-\nabla_Z\mathcal{L}_Vg-\nabla_{[V,Z]}g)(X,Y)&=&-g((\mathcal{L}_V\nabla)(X,Z),Y)\nonumber\\
      &&-g((\mathcal{L}_V\nabla)(Y,Z),X),\nonumber
    \end{eqnarray}
    where $g$ is the metric connection i.e., $\nabla g=0$. So the above equation reduces to,
    \begin{equation}\label{C6}
      (\nabla_Z\mathcal{L}_Vg)(X,Y)=g((\mathcal{L}_V\nabla)(X,Z),Y)+g((\mathcal{L}_V\nabla)(Y,Z),X),
    \end{equation}
    for all vector fields $X$, $Y$, $Z$ on $M$. Combining (\ref{C5}) and (\ref{C6}) and by a straightforward combinatorial computation and using the symmetry of $(\mathcal{L}_V\nabla)$ the foregoing equation yields,
    \begin{eqnarray}
      g((\mathcal{L}_V\nabla)(X,Y),Z) &=& (\nabla_ZS)(X,Y)-(\nabla_XS)(Y,Z)-(\nabla_YS)(Z,X) \nonumber\\
      && -2(\mu+1)\{g(X,Y)\eta(Z)-\eta(X)\eta(Y)\eta(Z)\},\label{C7}
    \end{eqnarray}
    for arbitrary vector fields $X$,$Y$ and $Z$ on $M$. Using (\ref{C1}) and (\ref{C2}), the foregoing equation yields,
    \begin{equation}\label{C8}
      (\mathcal{L}_V\nabla)(X,\xi)=2QX+4nX,
    \end{equation}
    for all $X\in\chi(M)$. Now differentiating covariantly this with respect to arbitrary vector field $Y$ we get,
    \begin{equation}\label{C9}
      (\nabla_Y\mathcal{L}_\nabla)(X,\xi)=2(\nabla_YQ)X-(\mathcal{L}_V)(X,Y)+\eta(Y)(2QX+4nX).
    \end{equation}
    We know that, $(\mathcal{L}_VR)(X,Y)Z=(\nabla_X\mathcal{L}_V\nabla)(Y,Z)- (\nabla_Y\mathcal{L}_V\nabla)(X,Z)$. Using (\ref{C9}) in the previous relation we get,
    \begin{eqnarray}\label{C10}
      (\mathcal{L}_VR)(X,Y)\xi&=&2\{(\nabla_XQ)Y-(\nabla_YQ)X\}+2\eta(X)(QY+2nY)\nonumber\\
      &&-2\eta(Y)(QX+2nX),
    \end{eqnarray}
    for arbitrary vector fields $X$ and $Y$ on $M$. Setting $Y=\xi$ in the aforementioned equation and using (\ref{B12}),(\ref{C1}) and (\ref{C2}) we get,
    \begin{equation}\label{C11}
      (\mathcal{L}_VR)(X,\xi)\xi=0.
    \end{equation}
    Now, taking (\ref{C4}) in account, the Lie derivative of $g(\xi,\xi)=1$ along the potential vector field $V$ gives rise to,
    \begin{equation}\label{C12}
      \eta(\mathcal{L}_V\xi)=\lambda-\frac{p}{2}-\frac{1}{(2n+1)}+\mu.
    \end{equation}
    Setting $Y=\xi$ and using (\ref{B2}) and (\ref{B6}) the equation (\ref{C4}) gives rise to,
    \begin{equation}\label{C13}
      (\mathcal{L}_V\eta)X-g(X,\mathcal{L}_V\xi)=(p+\frac{2}{(2n+1)}-2\lambda-2\mu)\eta(X),
    \end{equation}
    which holds for arbitrary vector field $X$ on $M$. From (\ref{B11}) we have, $R(X,\xi)\xi=\eta(X)\xi-X$. Taking Lie derivative along the potential vector field $V$ and taking (\ref{C12}) and (\ref{C13}) in account, this reduces to,
    \begin{equation}\label{C14}
      (\mathcal{L}_VR)(X,\xi)\xi=(2\lambda+2\mu-p-\frac{2}{(2n+1)})(X-\eta(X)\xi),
    \end{equation}
    for all $X\in\chi(M)$. Finally comparing (\ref{C11}) and (\ref{C14}) we have, $(2\lambda+2\mu-p-\frac{2}{(2n+1)})(X-\eta(X)\xi)=0$. Since this holds for arbitrary $X\in\chi(M)$ so, it reduces to,
    \begin{equation}\label{C15}
      \lambda=\frac{p}{2}+\frac{1}{(2n+1)}-\mu.
    \end{equation}
    Using the relation (\ref{C15}) in (\ref{C12}) we easily obtain, $\eta(\mathcal{L}_V\xi)=0$. Since we have considered the potential vector field $V$ as contact vector field so there must exists a smooth function $f$ such that $\mathcal{L}_V\xi=f\xi$. Making use of this in (\ref{C12}) we get $f=\lambda-\frac{p}{2}-\frac{1}{(2n+1)}+\mu$. Therefore by using the relation (\ref{C15}) we get $f=0$ and thus $\mathcal{L}_V\xi=0$. Finally the equation (\ref{C13}) reduces to,
    \begin{equation}\label{C16}
      \mathcal{L}_V\eta=0.
    \end{equation}
    So, $V$ is strictly infinitesimal contact transformation.\par
    We know the well-known formula from Yano\cite{Yano} that, $(\mathcal{L}_V\nabla)(X,Y)=\mathcal{L}_V\nabla_XY-\nabla_X\mathcal{L}_VY-\nabla_{[V,X]}Y$. Setting $Y=\xi$ and using (\ref{B9}), $\mathcal{L}_V\xi=0$ and (\ref{C16}) we get, $(\mathcal{L}_V\nabla)(X,\xi)=0$. Substituting this in (\ref{C8}) we finally obtain $QX=-2nX~\forall X\in\chi(M)$. This proves our result.
  \end{proof}
\end{thm}

$*$-conformal $\eta$-Ricci soliton is a mere generalisation of conformal $*$-Ricci soliton where we consider $\mu=0$ in (\ref{A1}) to get conformal $*$-Ricci soliton equation. We can rewrite the above theorem as:
\begin{cor}
Let $M^{(2n+1)}(\phi,\xi,\eta,g)$ be a Kenmotsu manifold. If the metric $g$ represents a conformal $*$-Ricci soliton and if the soliton vector field $V$ is contact, then V is strictly infinitesimal contact transformation and the manifold is Einstein.
\end{cor}

\begin{ex}
  Let us consider the set $M=\{(x,y,z,u,v)\in\Bbb{R}^5\}$ as our manifold where $(x,y,z,u,v)$ are the standard coordinates in $\Bbb{R}^5$. The vector fields defined below:
  \begin{align*}
    e_1 &= e^{-v}\frac{\partial}{\partial x}, & e_2 &= e^{-v}\frac{\partial}{\partial y}, & e_3 &= e^{-v}\frac{\partial}{\partial z}, & e_4 &= e^{-v}\frac{\partial}{\partial u}, & e_5 &= \frac{\partial}{\partial v}
  \end{align*}
  are linearly independent at each point of $M$. We define the metric $g$ as
  \[
  g(e_i,e_j)=\begin{cases}
               1, & \mbox{if } i=j~and~i,j\in\{1,2,5\} \\
               -1, & \mbox{if } i=j~and~i,j\in\{3,4\} \\
               0, & \mbox{otherwise}.
             \end{cases}
  \]
  Let $\eta$ be a 1-form defined by $\eta(X)=g(X,e_5)$, for arbitrary $X\in\chi(M)$. Let us define (1,1)-tensor field $\phi$ as:
  \begin{align*}
    \phi(e_1) &= e_3, & \phi(e_2) &= e_4, & \phi(e_3) &= -e_1, & \phi(e_4) &= -e_2, & \phi(e_5) &= 0.
  \end{align*}
    Then it satisfy the relations $\eta(\xi)=1$, $\phi^2(X)=-X+\eta(X)\xi$ and $g(\phi X,\phi Y)=g(X,Y)-\eta(X)\eta(Y)$ where $\xi=e_5$ and $X,Y$ is arbitrary vector field on $M$. So, $(M,\phi,\xi,\eta,g)$ defines an almost contact structure on $M$.\par
     We can now deduce that,
    \begin{align*}
      [e_1,e_2] &=0 & [e_1,e_3] &=0 & [e_1,e_4] &=0 & [e_1,e_5] &=e_1  \\
      [e_2,e_1] &=0 & [e_2,e_3] &=0 & [e_2,e_4] &=0 & [e_2,e_5] &=e_2  \\
      [e_3,e_1] &=0 & [e_3,e_2] &=0 & [e_3,e_4] &=0 & [e_3,e_5] &=e_3  \\
      [e_4,e_1] &=0 & [e_4,e_2] &=0 & [e_4,e_3] &=0 & [e_4,e_5] &=e_4  \\
      [e_5,e_1] &=-e_1 & [e_5,e_2] &=-e_2 & [e_5,e_3] &=-e_3 & [e_5,e_4] &=-e_4.
    \end{align*}
Let $\nabla$ be the Levi-Civita connection of $g$. Then from $Koszul's formula$ for arbitrary $X,Y,Z\in\chi(M)$ given by:\\
\begin{eqnarray}
2g(\nabla_XY,Z)&=&Xg(Y,Z)+Yg(Z,X)-Zg(X,Y)-g(X,[Y,Z])\nonumber\\
&-&g(Y,[X,Z])+g(Z,[X,Y]),\nonumber
\end{eqnarray}
we can have:
\begin{align*}
  \nabla_{e_1}e_1 &=-e_5 & \nabla_{e_1}e_2 &=0 & \nabla_{e_1}e_3 &=0 & \nabla_{e_1}e_4 &=0 & \nabla_{e_1}e_5 &=e_1 \\
  \nabla_{e_2}e_1 &=0 & \nabla_{e_2}e_2 &=-e_5 & \nabla_{e_2}e_3 &=0 & \nabla_{e_2}e_4 &=0 & \nabla_{e_2}e_5 &=e_2 \\
  \nabla_{e_3}e_1 &=0 & \nabla_{e_3}e_2 &=0 & \nabla_{e_3}e_3 &=-e_5 & \nabla_{e_3}e_4 &=0 & \nabla_{e_3}e_5 &=e_3 \\
  \nabla_{e_4}e_1 &=0 & \nabla_{e_4}e_2 &=0 & \nabla_{e_4}e_3 &=0 & \nabla_{e_4}e_4 &=-e_5 & \nabla_{e_4}e_5 &=e_4 \\
  \nabla_{e_5}e_1 &=0 & \nabla_{e_5}e_2 &=0 & \nabla_{e_5}e_3 &=0 & \nabla_{e_5}e_4 &=0 & \nabla_{e_5}e_5 &=0.
\end{align*}
Therefore $(\nabla_X\phi)Y=g(\phi X,Y)\xi-\eta(Y)\phi X$ is satisfied for arbitrary $X,Y\in\chi(M)$. So $(M,\phi,\xi,\eta,g)$ becomes a Kenmotsu manifold.\\
The non-vanishing components of curvature tensor are:
\begin{align*}
  R(e_1,e_2)e_2&=-e_1 & R(e_1,e_3)e_3&=-e_1 & R(e_1,e_4)e_4&=-e_1 \\
  R(e_1,e_5)e_5&=-e_1 & R(e_1,e_2)e_1&=e_2 & R(e_1,e_3)e_1&=e_3 \\
  R(e_1,e_4)e_1&=e_4 & R(e_1,e_5)e_1&=e_5 & R(e_2,e_3)e_2&=e_3 \\
  R(e_2,e_4)e_2&=e_4 & R(e_2,e_5)e_2&=e_5 & R(e_2,e_3)e_3&=-e_2 \\
  R(e_2,e_4)e_4&=-e_2 & R(e_2,e_5)e_5&=-e_2 & R(e_3,e_4)e_3&=e_4 \\
  R(e_3,e_5)e_3&=e_5 & R(e_3,e_4)e_4&=-e_3 & R(e_4,e_5)e_4&=e_5 \\ R(e_5,e_3)e_5&=e_3 & R(e_5,e_4)e_5&=e_4.
\end{align*}
Now from the above results we have, $S(e_i,e_i)=-4$ for $i=1,2,3,4,5$  and,
\begin{equation}\label{ex1 1}
  S(X,Y)=-4g(X,Y)~~\forall X,Y\in\chi(M).
\end{equation}
Contracting this we have $r=\sum_{i=1}^{5}S(e_i,e_i)=-20=-2n(2n+1)$ where dimension of the manifold $2n+1=5$. Also, we have,
\[
  S^*(e_i,e_i)=\begin{cases}
               -1, & \mbox{if } i=1,2,3,4 \\
               0, & \mbox{if } i=5.
               \end{cases}
\]
and, $r^*=r+4n^2=-20+16=-4$. So,
\begin{equation}\label{ex1 2}
  S^*(X,Y)=-g(X,Y)+\eta(X)\eta(Y)~~\forall X,Y\in\chi(M).
\end{equation}
Now we consider a vector field $V$ as,
\begin{equation}\label{ex1 3}
  V=x\frac{\partial}{\partial x}+y\frac{\partial}{\partial y}+z\frac{\partial}{\partial z}+u\frac{\partial}{\partial u}+\frac{\partial}{\partial v}.
\end{equation}
Then from the above results we can justify that,
\begin{equation}\label{ex1 4}
  (\mathcal{L}_Vg)(X,Y)=4\{g(X,Y)-\eta(X)\eta(Y)\},
\end{equation}
which holds for all $X,Y\in\chi(M)$. From (\ref{ex1 2}) and (\ref{ex1 4}) we can conclude that $g$ represents a $*$-conformal $\eta$-Ricci soliton i.e., it satisfies (\ref{A1}) for potential vector field $V$ defined by (\ref{ex1 3}), $\lambda=\frac{p}{2}-\frac{4}{5}$ and $\mu=1$.
\end{ex}

\begin{thm}
Let $M^{(2n+1)}(\phi,\xi,\eta,g)$ be a Kenmotsu manifold. If the metric $g$ represents a gradient almost $*$-conformal $\eta$-Ricci soliton then either $M$ is Einstein or there exists an open set where the potential vector field $V$ is pointwise collinear with the characteristic vector field $\xi$.
\begin{proof}
  Using (\ref{C3}) in the definition of gradient almost $*$-conformal $\eta$-Ricci soliton give by equation (\ref{A2}) we get,
  \begin{equation}\label{D1}
    \nabla_XDf=-QX-(\lambda-\frac{p}{2}-\frac{1}{(2n+1)}+2n-1)X-(\mu+1)\eta(X)\xi,
  \end{equation}
  for any vector field $X$ on $M$. Taking covariant derivative along arbitrary vector $Y$ and using (\ref{B9}), (\ref{B10}) we get,
  \begin{eqnarray}\label{D2}
    \nabla_Y\nabla_XDf &=& -(\nabla_YQ)X-Q(\nabla_YX)-Y(\lambda)X-(\lambda-\frac{p}{2}-\frac{1}{(2n+1)}\nonumber\\
    && +2n-1)(\nabla_YX)-(\mu+1)\{g(X,Y)\xi-2\eta(X)\eta(Y)\xi\nonumber\\
    &&+\eta(\nabla_YX)\xi+\eta(X)Y\}.
  \end{eqnarray}
  Applying this in the expression of Riemannian curvature tensor we get,
  \begin{eqnarray}\label{D3}
    R(X,Y)Df&=&(\nabla_YQ)X-(\nabla_XQ)Y+Y(\lambda)X-X(\lambda)Y\nonumber\\
    &&-(\mu+1)\{\eta(Y)X-\eta(X)Y\}.
  \end{eqnarray}
  Moreover an inner product w.r.t. $\xi$ and use of (\ref{C1}) and (\ref{C2}) yields,
  \begin{equation}\label{D4}
    g(R(X,Y)Df,\xi)=Y(\lambda)\eta(X)-X(\lambda)\eta(Y),
  \end{equation}
  for $X,Y\in\chi(M)$. Furthermore the inner product of (\ref{B11}) with the potential vector field $Df$ gives,
  \begin{equation}\label{D5}
    g(R(X,Y)Df,\xi)=\eta(Y)X(f)-\eta(X)Y(f),
  \end{equation}
  for arbitrary $X$ and $Y$ on $M$. Comparing (\ref{D4}) and (\ref{D5}) and setting $Y=\xi$ we have $X(f+\lambda)=\xi(f+\lambda)\eta(X)$. From this we obtain,
  \begin{equation}\label{D6}
    d(f+\lambda)=\xi(f+\lambda)\eta.
  \end{equation}
  So, $(f+\lambda)$ is invariant along the distribution $Ker(\eta)$ i.e., if $X\in Ker(\eta)$ then $X(f+\lambda)=d(f+\lambda)X=0$.\\
  Now, if we takhe inner product w.r.t. arbitrary vector field $Z$ after plugging $X=\xi$ in (\ref{D3}) we get,
  \begin{eqnarray}\label{D7}
    g(R(\xi,Y)Df,Z)&=&S(Y,Z)+(2n-\xi(\lambda)+\mu+1)g(Y,Z)+Y(\lambda)\eta(Z)\nonumber\\
    &&-(\mu+1)\eta(Y)\eta(Z).
  \end{eqnarray}
  Again from (\ref{B11}) we can easily obtain for arbitrary vector fields $Y$ and $Z$ on $M$,
  \begin{equation}\label{D8}
    g(R(\xi,Y)Df,Z)=\xi(f)g(Y,Z)-Y(f)\eta(Z).
  \end{equation}
  Comparing the equations (\ref{D7}) and (\ref{D8}) and using (\ref{D6}) we obtain,
  \begin{equation}\label{D9}
    S(Y,Z)=\{\xi(f+\lambda)-\mu-2n-1\}g(Y,Z)-\{\xi(f+\lambda)-\mu-1\}\eta(Y)\eta(Z).
  \end{equation}
  Since the above equation holds good for arbitrary $Y$ and $Z$, so the manifold is $\eta$-Einstein. Now contracting (\ref{D9}) we obtain,
  \begin{equation}\label{D10}
    \xi(f+\lambda)=\frac{r}{2n}+\mu+2n+2.
  \end{equation}
  Plugging this in (\ref{D9}) we acquire,
  \begin{equation}
    S(Y,Z)=(\frac{r}{2n}+1)g(Y,Z)-(\frac{r}{2n}+2n+1)\eta(Y)\eta(Z),\nonumber
  \end{equation}
  for arbitrary vector fields $Y$ and $Z$ on $M$ which is exactly same as (\ref{B15}). Now contracting (\ref{D3}) w.r.t. $X$ gives rise to,
  \begin{equation}\label{D11}
    S(Y,Df)=\frac{1}{2}Y(r)+2nY(\lambda)-2n(\mu+1)\eta(Y),
  \end{equation}
  which holds for any $Y\in\chi(M)$. Now, comparing this with (\ref{B15}) we obtain,
  \begin{eqnarray}
    && (r+2n)Y(f)-(r+2n(2n+1))\eta(Y)\xi(f)-nY(r) \nonumber\\
    && -4n^2Y(\lambda)+4n^2(\mu+1)\eta(y)=0,\label{D12}
  \end{eqnarray}
  for all $Y\in\chi(M)$. Now, setting $Y=\xi$ and then using (\ref{D10}) we easily obtain the relation,
  \begin{equation}\label{D13}
    \xi(r)=-2(r+2n(2n+1)).
  \end{equation}
  Since $d^2=0$ and $d\eta=0$, from (\ref{D6}) we obtain $dr\wedge\eta=0$ i.e., $dr(X)\eta(Y)-dr(Y)\eta(X)=0$ for arbitrary $X,Y\in\chi(M)$. After considering $Y=\xi$ and using (\ref{D13}) it reduces to $X(r)=-2(r+2n(2n+1))\xi$. Since $X$ is an arbitrary vector field so we conclude that,
  \begin{equation}\label{D14}
    Dr=-2(r+2n(2n+1))\xi.
  \end{equation}
  Let $X$ be a vector field of the distribution $Ker(\eta)$. Then, (\ref{D12}) reduces to,
  \begin{equation}\nonumber
    (r+2n)X(f)-4n^2X(\lambda)=0.
  \end{equation}
  Using (\ref{D6}) and (\ref{D10}) we obtain, $(r+2n(2n+1))X(f)=0$. From here we conclude,
  \begin{equation}\nonumber
    (r+2n(2n+1))(Df-\xi(f)\xi)=0.
  \end{equation}
  \par If $r=-2n(2n+1)$, then from (\ref{B15}) we get that the manifold is Einstein with Einstein constant $-2n$.\par
  If $r\neq-2n(2n+1)$ on some open set $O$ of $M$, then $Df=\xi(f)\xi$ on that open set that is, the potential vector field is pointwise collinear with the characteristic vector field $\xi$.

\end{proof}
\end{thm}

If we let the coefficient of $\eta\otimes\eta$ in (\ref{A2}) to be zero then the soliton reduces to gradient almost conformal $*$-Ricci soliton. The aforementioned result in the framework of gradient almost conformal $*$-Ricci soliton can be stated as:
\begin{cor}
Let $M^{(2n+1)}(\phi,\xi,\eta,g)$ be a Kenmotsu manifold. If the metric $g$ represents a gradient almost conformal $*$-Ricci soliton then either $M$ is Einstein or the potential vector field $V$ is pointwise collinear with the characteristic vector filed $\xi$ on an open set on $M$.

\end{cor}

\begin{ex}
  Let us consider the set $M=\{(x,y,z,u,v)\in\Bbb{R}^5\}$ as our manifold where $(x,y,z,u,v)$ are the standard coordinates in $\Bbb{R}^5$. The vector fields defined below:
  \begin{align*}
    e_1 &= v\frac{\partial}{\partial x}, & e_2 &= v\frac{\partial}{\partial y}, & e_3 &= v\frac{\partial}{\partial z}, & e_4 &= v\frac{\partial}{\partial u}, & e_5 &= -v\frac{\partial}{\partial v}
  \end{align*}
  forms a linearly independent set of vector fields on $M$. We define the metric $g$ as
  \begin{equation*}
    (g_{ij})=\begin{pmatrix}
               1 & 0 & 0 & 0 & 0 \\
               0 & 1 & 0 & 0 & 0 \\
               0 & 0 & 1 & 0 & 0 \\
               0 & 0 & 0 & 1 & 0 \\
               0 & 0 & 0 & 0 & 1
             \end{pmatrix}.
  \end{equation*}
   We consider the reeb vector field $\xi=e_5$ then the 1-form $\eta$ is defined by $\eta(X)=g(X,e_5)$, for arbitrary $X\in\chi(M)$ then, $\eta=dv$. Let us define (1,1)-tensor field $\phi$ as:
  \begin{align*}
    \phi(e_1) &= e_2, & \phi(e_2) &= -e_1, & \phi(e_3) &= e_4, & \phi(e_4) &= -e_3, & \phi(e_5) &= 0.
  \end{align*}
    Then it satisfy the relations $\eta(\xi)=1$, $\phi^2(X)=-X+\eta(X)\xi$ and $g(\phi X,\phi Y)=g(X,Y)-\eta(X)\eta(Y)$ where $X,Y$ is arbitrary vector field on $M$. So, $(M,\phi,\xi,\eta,g)$ defines an almost contact structure on $M$.\par
Let $\nabla$ be the Levi-Civita connection of $g$. Then from $Koszul's formula$ for arbitrary $X,Y,Z\in\chi(M)$ given by:\\
\begin{eqnarray}
2g(\nabla_XY,Z)&=&Xg(Y,Z)+Yg(Z,X)-Zg(X,Y)-g(X,[Y,Z])\nonumber\\
&-&g(Y,[X,Z])+g(Z,[X,Y]),\nonumber
\end{eqnarray}
we can have:
\begin{align*}
  \nabla_{e_1}e_1 &=-e_5 & \nabla_{e_1}e_2 &=0 & \nabla_{e_1}e_3 &=0 & \nabla_{e_1}e_4 &=0 & \nabla_{e_1}e_5 &=e_1 \\
  \nabla_{e_2}e_1 &=0 & \nabla_{e_2}e_2 &=-e_5 & \nabla_{e_2}e_3 &=0 & \nabla_{e_2}e_4 &=0 & \nabla_{e_2}e_5 &=e_2 \\
  \nabla_{e_3}e_1 &=0 & \nabla_{e_3}e_2 &=0 & \nabla_{e_3}e_3 &=-e_5 & \nabla_{e_3}e_4 &=0 & \nabla_{e_3}e_5 &=e_3 \\
  \nabla_{e_4}e_1 &=0 & \nabla_{e_4}e_2 &=0 & \nabla_{e_4}e_3 &=0 & \nabla_{e_4}e_4 &=-e_5 & \nabla_{e_4}e_5 &=e_4 \\
  \nabla_{e_5}e_1 &=0 & \nabla_{e_5}e_2 &=0 & \nabla_{e_5}e_3 &=0 & \nabla_{e_5}e_4 &=0 & \nabla_{e_5}e_5 &=0.
\end{align*}
Therefore $(\nabla_X\phi)Y=g(\phi X,Y)\xi-\eta(Y)\phi X$ is satisfied for arbitrary $X,Y\in\chi(M)$. So $(M,\phi,\xi,\eta,g)$ becomes a Kenmotsu manifold.\\
The non-vanishing components of curvature tensor are:
\begin{align*}
  R(e_1,e_2)e_2&=-e_1 & R(e_1,e_3)e_3&=-e_1 & R(e_1,e_4)e_4&=-e_1 \\
  R(e_1,e_5)e_5&=-e_1 & R(e_1,e_2)e_1&=e_2 & R(e_1,e_3)e_1&=e_3 \\
  R(e_1,e_4)e_1&=e_4 & R(e_1,e_5)e_1&=e_5 & R(e_2,e_3)e_2&=e_3 \\
  R(e_2,e_4)e_2&=e_4 & R(e_2,e_5)e_2&=e_5 & R(e_2,e_3)e_3&=-e_2 \\
  R(e_2,e_4)e_4&=-e_2 & R(e_2,e_5)e_5&=-e_2 & R(e_3,e_4)e_3&=e_4 \\
  R(e_3,e_5)e_3&=e_5 & R(e_3,e_4)e_4&=-e_3 & R(e_4,e_5)e_4&=e_5 \\ R(e_5,e_3)e_5&=e_3 & R(e_5,e_4)e_5&=e_4.
\end{align*}
Now from the above results we have, $S(e_i,e_i)=-4$ for $i=1,2,3,4,5$  and,
\begin{equation}\label{ex2 1}
  S(X,Y)=-4g(X,Y)~~\forall X,Y\in\chi(M).
\end{equation}
So, the manifold is Einstein. Also, we have,
\[
  S^*(e_i,e_i)=\begin{cases}
               -1, & \mbox{if } i=1,2,3,4 \\
               0, & \mbox{if } i=5.
               \end{cases}
\]
and,
\begin{equation}\label{ex2 2}
  S^*(X,Y)=-g(X,Y)+\eta(X)\eta(Y)~~\forall X,Y\in\chi(M).
\end{equation}
Let $f:M\rightarrow \mathbb{R}$ be a smooth function defined by,
\begin{equation}\label{ex2 3}
  f(x,y,z,u,v)=x^2+y^2+z^2+u^2+\frac{v^2}{2}.
\end{equation}
Then the gradient of $f$, $Df$ is given by,
\begin{equation}\label{ex2 4}
  Df=2x\frac{\partial}{\partial x}+2y\frac{\partial}{\partial y}+2z\frac{\partial}{\partial z}+2u\frac{\partial}{\partial u}+v\frac{\partial}{\partial v}.
\end{equation}
Then from the above results we can verify that,
\begin{equation}\label{ex2 5}
  (\mathcal{L}_{Df}g)(X,Y)=2\{g(X,Y)-\eta(X)\eta(Y)\},
\end{equation}
which holds for all $X,Y\in\chi(M)$. From (\ref{ex2 2}) and (\ref{ex2 5}) we obtain that $g$ represents a gradient almost $*$-conformal $\eta$-Ricci soliton i.e., it satisfies (\ref{A2}) for $V=Df$, where $f$ is defined by (\ref{ex2 3}), $\lambda=\frac{p}{2}+\frac{1}{5}$ and $\mu=0$.
\end{ex}

\section{\textbf{$*$-conformal $\eta$-Ricci soliton on $(\kappa,\mu)^\prime$-almost Kenmotsu manifold with $\kappa<-1$}}
In this section we consider the manifold as a $(2n+1)$-dimensional almost Kenmotsu manifold where the characteristic vector field $\xi$ satisfies $(\kappa,-2)^\prime$-nullity distribution. Then we let the metric $g$ to represent a $*$-conformal $\eta$-Ricci soliton. Here we look back on some pertinent results and used these in our work.

\begin{lem}\cite{Dai}
On a $(\kappa,-2)^\prime$-almost Kenmotsu manifold with $\kappa<-1$ the $*$-Ricci tensor is given by
\begin{equation}\label{E1}
  S^*(X,Y)=-(\kappa+2)(g(X,Y)-\eta(X)\eta(Y)),
\end{equation}
for any vector fields $X$ and $Y$.
\end{lem}

\begin{thm}
Let $M^{(2n+1)}(\phi,\xi,\eta,g)$ be an almost Kenmotsu manifold such that $\xi$ belongs to $(\kappa,-2)^\prime$-nullity distribution where $\kappa<-1$. If the metric $g$ represents a $*$-conformal $\eta$-Ricci soliton satisfying $p\neq2\lambda+2\mu-\frac{2}{(2n+1)}$ then, $M$ is Ricci-flat and is locally isometric to $\mathbb{H}^{n+1}(-4)\times\mathbb{R}^n$.
\begin{proof}
  Combining (\ref{A1}) with (\ref{E1}) we get,
  \begin{equation}\label{E2}
    (\mathcal{L}_Vg)(X,Y)=(p+2\kappa-2\lambda+4+\frac{2}{(2n+1)})g(X,Y)-2(\kappa+\mu+2)\eta(X)\eta(Y),
  \end{equation}
  for all vector fields $X$ and $Y$ on $M$. Now taking covariant derivative of the foregoing equation along arbitrary vector field $Z$ and using (\ref{B25}) we get,
  \begin{eqnarray}\nonumber
    (\nabla_Z\mathcal{L}_Vg)(X,Y)&=&-2(\kappa+\mu+2)[\eta(Y)g(X,Z)+\eta(X)g(Y,Z)+\eta(Y)\nonumber\\
    &&g(h^\prime Z,X)+\eta(X)g(h^\prime Z,Y)-2\eta(X)\eta(Y)\eta(Z)].
  \end{eqnarray}
   By a straightforward combinatorial computation, use of (\ref{C6}), the symmetry of $(\mathcal{L}_V\nabla)$ in the aforementioned equation we get,
   \begin{equation}\label{E3}
     (\mathcal{L}_V\nabla)(X,Y)=-2(\kappa+\mu+2)[g(X,Y)+g(h^\prime X,Y)-\eta(X)\eta(Y)]\xi,
   \end{equation}
   for all $X,Y\in\chi(M)$. Replacing $Y=\xi$ and using (\ref{B2}), (\ref{B6}) and (\ref{B16}) we have,
   \begin{equation}\label{E4}
     (\mathcal{L}_V\nabla)(X,\xi)=0,
   \end{equation}
   for arbitrary vector field $X$ on $M$. Now taking (\ref{B18}) and (\ref{E3}) into account and differentiating (\ref{E4}) covariantly along arbitrary vector field $Y$ we get,
   \begin{equation}\label{E5}
     (\nabla_Y\mathcal{L}_V\nabla)(X,\xi)=2(\kappa+\mu+2)[g(X,Y)-\eta(X)\eta(Y)+2g(h^\prime X,Y)+g(h^{\prime 2}X,Y)]\xi
   \end{equation}
   for any vector fields $X$ and $Y$ on $M$. Again from Yano we have the well-known curvature property, $(\mathcal{L}_VR)(X,Y)Z=(\nabla_X\mathcal{L}_V\nabla)(Y,Z)-(\nabla_Y\mathcal{L}_V\nabla)(X,Z)$. Setting $Z=\xi$ and using (\ref{E5}) repeatedly we obtain,
   \begin{equation}\label{E6}
     (\mathcal{L}_VR)(X,Y)\xi=0,
   \end{equation}
   for arbitrary $X,Y\in\chi(M)$. Now taking lie derivative of (\ref{B20}) along the potential vector field $V$, taking (\ref{B2}) and (\ref{B16}) into account we get,
   \begin{eqnarray}\label{E7}
     (\mathcal{L}_VR)(X,\xi)\xi &=& \kappa[g(X,\mathcal{L}_V\xi)\xi-2\eta(\mathcal{L}_V\xi)X-((\mathcal{L}_V\eta)X)\xi]+2[2\eta(\mathcal{L}_V\xi) \nonumber \\
     && h^\prime X-\eta(X)(h^\prime(\mathcal{L}_V\xi))-g(h^\prime X,\mathcal{L}_V\xi)\xi-((\mathcal{L}_Vh^\prime)X)],\nonumber\\
     &&
   \end{eqnarray}
   for any vector field $X$ on $M$. Plugging $Y=\xi$ in (\ref{E2}) we obtain,
   \begin{equation}\label{E8}
     (\mathcal{L}_V\eta)X-g(X,\mathcal{L}_V\xi)=(p-2\lambda-2\mu+\frac{2}{(2n+1)})\eta(X),
   \end{equation}
   for all $X\in\chi(M)$. Setting $X=\xi$ in the foregoing equation we get,
   \begin{equation}\label{E9}
     \eta(\mathcal{L}_V\xi)=-(\frac{p}{2}-\lambda-\mu+\frac{1}{(2n+1)}).
   \end{equation}
   With the help of (\ref{E6}),(\ref{E8}) and (\ref{E9}) we can rewrite the equation (\ref{E7}) as,
   \begin{eqnarray}\label{E10}
     && \kappa(p-2\lambda-2\mu+\frac{2}{(2n+1)})(X-\eta(X)\xi)-2(p-2\lambda-2\mu+\frac{2}{(2n+1)})h^\prime X-\nonumber \\
     &&2\eta(X)h^\prime(\mathcal{L}_V\xi)-2g(h^\prime X,\mathcal{L}_V\xi)\xi-2(\mathcal{L}_Vh^\prime)X=0.
   \end{eqnarray}
   Taking inner product of the foregoing equation with arbitrary vector field $Y$ on $M$ we obtain,
   \begin{gather}
     (p-2\lambda-2\mu+\frac{2}{(2n+1)})[\kappa(g(X,Y)-\eta(X)\eta(Y))-2g(h^\prime X,Y)]\nonumber \\
     -2\eta(X)g(h^\prime(\mathcal{L}_V\xi),Y)-2g(h^\prime X,\mathcal{L}_V\xi)\eta(Y)-2g((\mathcal{L}_Vh^\prime)X,Y)=0.\nonumber\\
   \end{gather}
   Since the above equation holds for any vector fields $X$ and $Y$ on $M$, by replacing $X$ by $\phi(X)$ and $Y$ by $\phi(Y)$ and taking (\ref{B5}) into account we get,
   \begin{equation}\label{E11}
     (p-2\lambda-2\mu+\frac{2}{(2n+1)})[\kappa g(\phi X,\phi Y)-2g(h^\prime \phi X,\phi Y)]-2g((\mathcal{L}_Vh^\prime)\phi X,\phi Y)=0,
   \end{equation}
   for all $X,Y\in\chi(M)$. Since $spec(h^\prime)=\{0,\alpha,-\alpha\}$, let $X$ and $V$ belong to the eigenspaces of $-\alpha$ and $\alpha$ denoted by $[-\alpha]^\prime$ and $[\alpha]^\prime$ respectively. Then $\phi X\in[\alpha]^\prime$ (for more details we refer to \cite{Dileo}). Then (\ref{E11}) can be rewritten as,
   \begin{equation}\label{E12}
     (p-2\lambda-2\mu+\frac{2}{(2n+1)})(\kappa-2) g(\phi X,\phi Y)-2g((\mathcal{L}_Vh^\prime)\phi X,\phi Y)=0,
   \end{equation}
   for all $X,Y\in\chi(M)$. It is remained to find the value of $g((\mathcal{L}_Vh^\prime)\phi X,\phi Y)$. To get this we prove a more generalized result: In a $(\kappa,\mu)^\prime$-almost Kenmotsu manifold $(\mathcal{L}_Xh^\prime)Y=0$, where $X$ and $Y$ belong to same eigenspaces.\par
   Without loss of generality we assume that $X,Y\in[\alpha]^\prime$ where $spec(h^\prime)=\{0,\alpha,-\alpha\}$. If we consider a local orthonormal $\phi$-basis as $\{\xi,e_i,\phi e_i\},i=1,2,...,n$ then,
   \begin{equation}\nonumber
     \nabla_XY=\sum_{i=1}^{n}g(\nabla_XY,e_i)e_i-(\alpha+1)g(X,Y)\xi.
   \end{equation}
   and,
   \begin{eqnarray}
     (\mathcal{L}_Xh^\prime)Y &=& \mathcal{L}_X(h^\prime Y)-h^\prime (\mathcal{L}_XY)\nonumber \\
     &=& \alpha (\mathcal{L}_XY)-h^\prime(\mathcal{L}_XY)\nonumber\\
     &=& \alpha (\nabla_XY-\nabla_YX)-h^\prime(\nabla_XY-\nabla_YX)\nonumber\\
     &=& \alpha(\alpha+1)g(X,Y)\xi-\alpha(\alpha+1)g(X,Y)\xi\nonumber\\
     &=& 0.\nonumber
   \end{eqnarray}
   Similarly we can prove that the above results hold good if $X,Y\in[-\alpha]^\prime$. For more details we refer to \cite{Dileo}. Now (\ref{E12}) reduces to,
   \begin{equation}\label{E13}
     (p-2\lambda-2\mu+\frac{2}{(2n+1)})(\kappa-2) g(\phi X,\phi Y)=0,
   \end{equation}
   for any vector fields $X$ and $Y$ on $M$. Since by hypothesis $p\neq2\lambda+2\mu-\frac{2}{(2n+1)}$, from the foregoing equation we infer that $\kappa=2\alpha$. Again form $\alpha^2=-(\kappa+1)$ we get $\alpha=-1$ and $\kappa=-2$. Plugging the value of $\kappa$ in (\ref{E1}) we have $S^*=0$, i.e., the manifold is Ricci-flat.\par
    Again we get $spec(h^\prime)=\{0,1,-1\}$. From corollary 4.2 of \cite{Dileo} we get $M$ is locally symmetric. From proposition 4.1 of \cite{Dileo} we finally conclude that $M$ is locally isometric to $\mathbb{H}^{n+1}(-4)\times\mathbb{R}^n$, where $\mathbb{H}^{n+1}(-4)$ is the hyperbolic space of constant curvature $-4$.
\end{proof}
\end{thm}

As we know, setting $\mu=0$ in (\ref{A1}) gives rise to the equation of conformal $*$-Ricci soliton, we can revisit the theorem-4.2 and can note the statement as:
\begin{cor}
Let $M(\phi,\xi,\eta,g)$ be a $(2n+1)$-dimensional almost Kenmotsu manifold such that $\xi$ blongs to $(\kappa,-2)^\prime$-nullity distribution where $\kappa<-1$. If the metric $g$ represents a conformal $*$-Ricci soliton satisfying $p\neq2\lambda-\frac{2}{(2n+1)}$ then, $M$ is Ricci-flat and is locally isometric to $\mathbb{H}^{n+1}(-4)\times\mathbb{R}^n$.
\end{cor}

\begin{ex}
We consider the manifold as $M=\{(x,y,z)\in\mathbb{R}^3:y\neg0\}$. We define three vector fields $e_1$, $e_2$ and $e_3$ as:
  \begin{align*}
    e_1 &= \frac{\partial}{\partial x}, & e_2 &= \frac{\partial}{\partial y}, & e_3 &= 2x\frac{\partial}{\partial x}-\frac{\partial}{\partial y}+\frac{\partial}{\partial z}.
  \end{align*}
  Then the set $\{e_1,e_2,e_3\}$ forms a linearly independent set of vector fields on $M$. We define the metric $g$ as
  \begin{equation*}
    (g_{ij})=\delta_{ij}~~\forall i,j\in\{1,2,3\}.
  \end{equation*}
   Then it is easy to verify that $\{e_1,e_2,e_3\}$ forms an orthonormal basis on $M$. Let the 1-form $\eta$ be defined by $\eta(X)=g(X,e_3)$, for arbitrary $X\in\chi(M)$. Let us define (1,1)-tensor field $\phi$ as:
  \begin{align*}
    \phi(e_1) &= -e_2, & \phi(e_2) &= e_1, & \phi(e_3) &= 0.
  \end{align*}
    Then it satisfy the relations $\eta(\xi)=1$, $\phi^2(X)=-X+\eta(X)\xi$ and $g(\phi X,\phi Y)=g(X,Y)-\eta(X)\eta(Y)$ where $\xi=e_3$ and $X,Y$ is arbitrary vector field on $M$. So, $(M,\phi,\xi,\eta,g)$ defines an almost contact structure on $M$.\par
We also can compute that,
    \begin{align*}
      [e_1,e_2] &=0, & [e_2,e_3] &=0, & [e_1,e_3] &=2e_1.
    \end{align*}
Let $\nabla$ be the Levi-Civita connection of $g$. Then from $Koszul's formula$ for arbitrary $X,Y,Z\in\chi(M)$ given by:\\
\begin{eqnarray}
2g(\nabla_XY,Z)&=&Xg(Y,Z)+Yg(Z,X)-Zg(X,Y)-g(X,[Y,Z])\nonumber\\
&-&g(Y,[X,Z])+g(Z,[X,Y]),\nonumber
\end{eqnarray}
we can have:
\begin{align*}
  \nabla_{e_1}e_1 &=-2e_3 & \nabla_{e_1}e_2 &=0 & \nabla_{e_1}e_3 &=2e_1 \\
  \nabla_{e_2}e_1 &=0 & \nabla_{e_2}e_2 &=0 & \nabla_{e_2}e_3 &=0 \\
  \nabla_{e_3}e_1 &=0 & \nabla_{e_3}e_2 &=0 & \nabla_{e_3}e_3 &=0.
\end{align*}
Therefore it is easy to verify that the structure $(M,\phi,\xi,\eta,g)$ is not Kenmotsu manifold. Now let us define the operator $h^\prime$ as:
    \begin{align*}
      h^\prime(e_1) &=e_1, & h^\prime(e_2) &=-e_2, & h^\prime(e_3) &=0.
    \end{align*}
By straightforward computation we have the components of curvature tensor as:
\begin{align*}
  R(e_1,e_2)e_1&=0 & R(e_1,e_2)e_2&=0 & R(e_1,e_2)e_3&=0 \\
  R(e_2,e_3)e_1&=0 & R(e_2,e_3)e_2&=0 & R(e_2,e_3)e_3&=0 \\
  R(e_1,e_3)e_1&=4e_3 & R(e_1,e_3)e_2&=0 & R(e_1,e_3)e_3&=-4e_1.
\end{align*}
Now from the above results and taking (\ref{B20}) in account we conclude that the reeb vector field $\xi$ belongs to the $(\kappa,-2)^\prime$-nullity distribution with $\kappa=-2$. So, the manifold is $(-2,-2)^\prime$-almost Kenmotsu manifold.\par
Now from (\ref{E1}) we get $S^*(X,Y)=0~~\forall X,Y\in\chi(M)$.
Now we consider a vector field $V$ as,
\begin{equation}\label{ex3 1}
  V=e^{2z}\frac{\partial}{\partial x}+4(y+z)\frac{\partial}{\partial y}.
\end{equation}
Then from the above results one can justify that,
\begin{align*}
  (\mathcal{L}_Vg)(e_1,e_1)&=0 & (\mathcal{L}_Vg)(e_2,e_2)&=8 & (\mathcal{L}_Vg)(e_3,e_3)&=0 \\
  (\mathcal{L}_Vg)(e_1,e_2)&=0 & (\mathcal{L}_Vg)(e_2,e_3)&=0 & (\mathcal{L}_Vg)(e_1,e_3)&=0.
\end{align*}
From here we can conclude that $g$ represents a $*$-conformal $\eta$-Ricci soliton i.e., it satisfies (\ref{A1}) for potential vector field $V$ defined by (\ref{ex3 1}), $\lambda=\frac{p}{2}-\frac{11}{3}$ and $\mu=4$. From Dileo and Pastore\cite{Dileo} we can further conclude that the manifold is locally isometric to $\mathbb{H}^2(-4)\times\mathbb{R}$.
\end{ex}

\section{\textbf{Conclusion}}
In this article, we have used methods of local Riemannian geometry to interpretation solutions of (\ref{A2}) and permeate Einstein metrics in a huge class of metrics as $*$-conformal $\eta$-Ricci solitons and almost $*$-conformal $\eta$-Ricci solitons on contact geometry, specially on Kenmotsu manifold. Not only our results will play important and motivational role in contact geometry but also there are further scope of research in this direction within the framework of various complex manifolds like K$\ddot{a}$hler manifold etc. There are some questions arises from our article to study in further research:\par
(i) Is theorem 3.3 true without assuming that the soliton vector field is contact?\par
(ii) Are the results of theorem 4.2 true if the characteristic vector field $\xi$ does not belong to $(\kappa,\mu)^\prime$-nullity distribution or if $\xi$ belongs to $(\kappa,\mu)^\prime$-nullity distribution with $\kappa\geq-1$?\par
(iii) Which of the results of our paper are also true for nearly Kenmotsu manifolds or $f$-Kenmotsu manifolds?

\section{\textbf{Acknowledgements}}
The authors are very much thankful to Vladimir Rovenski for his valuable comments and suggestions for improvements of this paper and for his constant support throughout the preparation of this paper. The first author is the corresponding author and this work was financially supported by UGC Senior Research Fellowship of India, Sr. No. 2061540940. Ref. No:21/06/2015(i)EU-V.

\end{document}